\documentclass[10pt]{amsart}

\usepackage{amsmath}
\usepackage{amsthm}
\usepackage{amsfonts}
\usepackage{amssymb,latexsym}
\usepackage[all]{xy}
\usepackage{graphicx}
\usepackage[mathscr]{eucal}
\usepackage{verbatim}
\usepackage{hyperref} 

\renewcommand{\baselinestretch}{1.08}

\theoremstyle{plain}
    \newtheorem{thm}{Theorem}[section]
    \newtheorem{prop}[thm]{Proposition}

    \newtheorem{cor}[thm]{Corollary}

\theoremstyle{definition}

\theoremstyle{remark}
    \newtheorem{rem}[thm]{Remark}
    \newtheorem{example}[thm]{Example}

\numberwithin{equation}{section}

\renewcommand{\Im}{\textup{Im}}

\newcommand{\ints}{\ensuremath{\mathbb{Z}}}
\newcommand{\rar}{\ensuremath{\rightarrow}}
\newcommand{\lrar}{\ensuremath{\longrightarrow}}

\newcommand{\la}{\langle}
\newcommand{\ra}{\rangle}
\newcommand{\Hom}{\textup{Hom}}
\newcommand{\StMod}{\textup{StMod}}
\newcommand{\stmod}{\textup{stmod}}
\newcommand{\Mod}{\textup{Mod}}

\newcommand{\uHom}{\underline{\Hom}}
\newcommand{\Ker}{\textup{Ker}}

\newcommand{\mcP}{\mathcal P}
\newcommand{\G}{\mathcal G}
\newcommand{\T}{\mathcal{T}}
\newcommand{\stk}[1]{\stackrel{#1}{\rightarrow}}
\newcommand{\lstk}[1]{\stackrel{#1}{\longrightarrow}}

\newcommand{\Tate}{\widehat{H}}

\newcommand{\Ann}{\textup{Ann}}
\newcommand{\ulp}{\textup{(}}
\newcommand{\urp}{\textup{)}}

\newcommand{\wOmega}{\widetilde{\Omega}}

\newcommand{\Ind}{\textup{Ind}}

\newcommand{\uInd}{\underline{\Ind}}

\newcommand{\up}{{\uparrow}^G}
\newcommand{\down}{{\downarrow}_H}

\newcommand{\Soc}{\textup{Soc}}
\newcommand{\Rad}{\textup{Rad}}

\begin{document}

\title{Ghosts in Modular representation theory}
\date{\today}

\author{Sunil K. Chebolu}
\address{Department of Mathematics \\
University of Western Ontario \\
London, ON N6A 5B7, Canada}
\email{schebolu@uwo.ca}

\author{J. Daniel Christensen}
\address{Department of Mathematics \\
University of Western Ontario \\
London, ON N6A 5B7, Canada}
\email{jdc@uwo.ca}

\author{J\'{a}n Min\'{a}\v{c}}
\address{Department of Mathematics\\
University of Western Ontario\\
London, ON N6A 5B7, Canada}
\email{minac@uwo.ca}


\begin{abstract}
A \emph{ghost} over a finite $p$-group $G$ is a map between modular
representations of $G$ which is invisible in Tate cohomology. Motivated by the
failure of the \emph{generating hypothesis}---the statement that ghosts
between finite-dimensional $G$-representations factor through a
projective---we define the \emph{ghost number} of $kG$ to be the smallest integer $l$
such that the composite of any $l$ ghosts between finite-dimensional
$G$-representations factors through a projective. In this paper we study ghosts
and the ghost numbers of $p$-groups. We begin by showing that a weaker version
of the generating hypothesis, where the target of the ghost is fixed to be the
trivial representation $k$, holds for all $p$-groups.
We then compute the ghost
numbers of all cyclic $p$-groups and all abelian $2$-groups with $C_2$ as a
summand. We obtain bounds on the ghost numbers for abelian $p$-groups and for
all $2$-groups which have a cyclic subgroup of index $2$. Using these bounds
we determine the finite abelian groups which have ghost number at most $2$.
Our methods involve techniques from group theory, representation theory,
triangulated category theory, and constructions motivated from homotopy
theory.
\end{abstract}

\vspace*{8pt}
\maketitle
\thispagestyle{empty}


\renewcommand{\baselinestretch}{1.08}

\section{Introduction}
Let $G$ be a finite $p$-group and let $k$ be a field of characteristic $p$.
Recall that the \emph{stable module category} $\StMod(kG)$ is the following
tensor triangulated category.  The objects are
left $kG$-modules and the space of morphisms between $kG$-modules $M$ and
$N$, denoted $\uHom_{kG}(M,N)$, is the $k$-vector space of $kG$-module
homomorphisms modulo those maps that factor through a projective module. The
category $\stmod(kG)$ is the full subcategory of finite-dimensional left $kG$-modules.
A \emph{ghost} in the stable module category 
is a map between $kG$-modules that is trivial in Tate cohomology.  In~\cite{CCM3}, 
we formulated the \emph{generating hypothesis} (GH) for $kG$ as
the statement that all ghosts between finite-dimensional $kG$-modules are
trivial in the stable module category, i.e., they factor through a projective.
(This formulation was motivated by the famous classical generating hypothesis
of Peter Freyd \cite{freydGH} in the stable homotopy category which is the conjecture that
there no non-trivial maps between finite spectra that are trivial in stable
homotopy groups.) We have shown in~\cite{CCM3} that the GH holds for $kG$, where
$G$ is a non-trivial finite $p$-group and $k$ is a field of characteristic $p$,
if and only if $G$ is either $C_2$ or $C_3$.

Motivated by the failure of the GH, we proceed in two natural directions. The first
one addresses the GH with the trivial representation $k$ as the
target. More precisely, we show that in the stable module category of any
$p$-group, a map $M \rar k$ from a finite-dimensional module to the trivial
module is stably trivial if it induces the zero map in Tate cohomology. We
give two proofs of this result in Section~\ref{sec:GH}. One of them uses
Spanier-Whitehead duality and Tate duality to
show that a map between finite-dimensional modules over a $p$-group is a
ghost if and only if it induces
the zero map on the functor $\uHom_{kG}(-, \Omega^* k )$.

The second direction we take measures the degree to which the
GH fails in $p$-groups other than $C_2$ and $C_3$. 
We define the \emph{ghost number} of $kG$ to be
the smallest non-negative integer $l$ such that the composite of any $l$ ghosts between
finite-dimensional $kG$-modules is trivial.  This is a new invariant for a group algebra. (Note that, in this terminology,
$C_2$ and $C_3$ are the only $p$-groups whose group algebras have ghost number $1$.)
%
%
A concept that will be key to our analysis of ghost numbers and related
invariants is that of a \emph{projective class}; see Section~\ref{sec:projectiveclass} 
for the definition.  
In Section~\ref{sec:universalghosts} we prove the existence of universal ghosts 
which implies that ghosts form part of a projective class.
This allows us to derive bounds on the ghost number, such
as Theorem~\ref{thm:3inequalitiesofnumbers} which states that the
ghost number of $kG$ is strictly less than the nilpotency
index of the Jacobson radical $J(kG)$ of $kG$.
Using these bounds, in Section~\ref{sec:computing}
we compute the ghost numbers of some $p$-groups. We show that the
ghost number of $kC_{p^r}$ is $\lceil(p^r - 1)/2 \rceil$, where
$\lceil x \rceil$ is the smallest integer that is greater than or equal to
$x$. If $G$ is a finite abelian $2$-group with $C_2$ as a summand, then the
ghost number of $kG$ is shown to be one less than the nilpotency index
of $J(kG)$. Computing the ghost number of
an arbitrary group algebra seems to be a hard problem.
Experience tells us that finding lower bounds for ghost numbers is much harder than
finding upper bounds. We obtain reasonable lower bounds for the ghost
numbers of the group algebras of abelian $p$-groups. We use these bounds to
show that the only abelian $p$-groups with ghost number $2$ are $C_4$,
$C_2 \oplus C_2$, and $C_5$.


The proofs of the aforementioned results involve a pleasant mix of methods
from group theory, representation theory and triangulated category theory.

Similar results on \emph{phantom maps} (maps between $kG$-modules which
factor through a projective when restricted to finite-dimensional
modules) in the stable module category appear in the work of Benson and
Gnacadja; see~\cite{bengna-1}.



We assume throughout that the group $G$ is a finite $p$-group.
The characteristic of the field $k$ is always assumed to divide the order of the group. For example,
when we write $kC_3$, the reader will understand that the characteristic of
$k$ is $3$.  When we speak of suspensions of a $kG$-module $M$, we mean
$\Omega^i\,M$ for any integer $i$, as an object of $\StMod(kG)$. We write
$\wOmega^i\,M$ for the projective-free part of $\Omega^i\,M$, a well-defined
$kG$-module. When we speak of Heller shifts of $M$, we mean $\wOmega^i\,M$ for
any integer $i$. We use standard facts about the stable module category and $kG$-modules; a good reference
is \cite{carlson-modulesandgroupalgebras}.

We would like to thank Dave Benson~\cite{CCM3} for 
Propositions~\ref{prop:conjugacyghost} and~\ref{prop:GHranktwogroups}
which allowed us to strengthen some of our results.
We thank Mark Hovey and Keir Lockridge for helpful conversations
and an anonymous referee for comments used to improve the exposition.

\section{The generating hypothesis} \label{sec:GH}

A map $\phi\colon M \rar N$ between $kG$-modules is said to be a \emph{ghost} if
the induced map
\[ \uHom_{kG}(\Omega^i k, M) \lrar \uHom_{kG}(\Omega^i k, N)\]
between the Tate cohomology groups is zero for each integer $i$.
(Recall that the Tate group $\Tate^i(G,M)$ of $G$ with coefficients in $M$ is
isomorphic to $\uHom_{kG}(\Omega^i k, M)$.)
If $G$ is a finite $p$-group and $k$ is a field
of characteristic $p$, then ``the generating hypothesis for $kG$'' is the
statement that all ghosts between finite-dimensional $kG$-modules are trivial
in the stable module category. In \cite{CCM3} we have shown that the only
non-trivial  $p$-groups for which the GH holds are the cyclic groups $C_2$ and
$C_3$.

\subsection{The generating hypothesis with target $k$}
While the generating hypothesis generally fails in the stable module
category, we show that a weak version of it holds for all $p$-groups. We begin
with motivation coming from homotopy theory for studying this weak
version.  Devinatz~\cite{ethanGH} proved the following partial affirmative
result on the generating hypothesis for the $p$-local stable homotopy category
of spectra where $p$ is an odd prime: If $\phi \colon X \rar S^0$ is a map from a
 finite spectrum to the sphere spectrum such that $\pi_* (\phi) = 0$, then
the $K_{(p)}$-localisation of $\phi$
 is trivial, where $K_{(p)}$ is periodic complex $K$-theory localised at $p$.

Motivated by this result, we consider  ``the GH with target $k$'',
which is the
statement that every map $M \rar k$ from a finite-dimensional $kG$-module $M$
to the trivial representation $k$ that induces the zero map in Tate cohomology
is trivial in the stable module category. We show that the GH with target $k$
holds for all $p$-groups. In fact, we give two proofs of this fact below; see
Corollaries~\ref{cor:GHwithtargetk-1} and~\ref{cor:GHwithtargetk-2}.




\begin{prop} \label{prop:H0H-1}
Let $G$ be a $p$-group and let $f \colon U \rar V$ be a ghost between
projective-free $kG$-modules.   Then we have the following:
\begin{enumerate}
 \item $U^G$ is contained in $\Ker(f)$. ($U^G$ is the $G$-invariant submodule of $U$.)
\item $\Im(f)$ is contained in $JV$. ($J$ denotes the Jacobson radical of $kG$.)
\end{enumerate}
\end{prop}

\begin{proof}
Since $f$ is a ghost, the induced map in  Tate cohomology is zero. In
particular, the maps
\begin{equation*}
 \Tate^0(G, f) \colon \frac{U^G}{\Im(N)} \lrar \frac{V^G}{\Im(N)}
\end{equation*}
and
\begin{equation*}
 \Tate^{-1}(G, f) \colon \frac{\Ker(N)}{JU} \lrar
\frac{\Ker(N)}{JV}
\end{equation*}
are zero. (Here $\Im(N)$ and $\Ker(N)$ are respectively  the image and kernel
of norm maps.) Since the norm map is trivial on  any projective-free module,
the above maps can be written as
\begin{equation}
\label{eqn:H0} \Tate^0(G, f) \colon {U^G} \lrar {V^G}
\end{equation}
and
\begin{equation}
\label{eqn:H-1} \Tate^{-1}(G, f) \colon \frac{U}{JU} \lrar \frac{V}{JV}.
\end{equation}
Both of the above maps are zero. The first part of the
proposition follows from (\ref{eqn:H0}) and the second part from
(\ref{eqn:H-1}).
\end{proof}

\begin{cor} \label{cor:GHwithtargetk-1}
The GH with target $k$ holds for $p$-groups.
\end{cor}

\begin{proof}
Let $f \colon M \rar k$ be a ghost in $\stmod(kG)$. $M$ is isomorphic in $\stmod(kG)$ to a projective-free
$kG$-module. Therefore we may assume that
$M$ is projective-free. Then, by Proposition \ref{prop:H0H-1}, the image of  $f \colon M \rar
k$ is contained in $J(k)$, which is zero. So we are done.
\end{proof}

Recall that associated to a $kG$-module $M$, one
has the socle (ascending) series $\Soc^i M$ and the  radical (descending)
series $\Rad^i M$.  For $i = 1$, $\Soc^1 M := M^G$, the $G$-invariant
submodule of $M$, and for $i > 1$, $\Soc^i M$ is defined inductively by
$\Soc^i M/\Soc^{i -1} M \cong (M / \Soc^{i -1} M)^G$.  $\Rad^i M := J^i M$ for
all $i \ge 0$, where $J^i$ denotes the $i$th power of the Jacobson radical of
$kG$. See~\cite{ben-1} for some properties of these series.
Proposition~\ref{prop:H0H-1} can now  be generalised as follows.

\begin{cor}\label{cor:soc-rad}
Let $G$ be a $p$-group and let $f \colon U \rar V$ be a composite of 
$l$ ghosts between projective-free $kG$-modules.
Then we have the following:
\begin{enumerate}
\item $\Soc^l(U)$ is contained in $\Ker(f)$.
\item $\Im(f)$ is contained in $\Rad^l(V)$.
\end{enumerate}
\end{cor}

\begin{proof}
This follows by a straightforward induction using
Proposition~\ref{prop:H0H-1}.
\end{proof}

\subsection{Ghosts and  duality}

A map $d \colon M \rar N$ between $kG$-modules is
called a \emph{dual ghost} if the induced map
\[\uHom_{kG}(M,  \Omega^i k) \longleftarrow \uHom_{kG}(N , \Omega^i k)\]
is zero for all $i$. Recall that for every  $kG$-module $L$, there is a
corresponding dual $kG$-module $L^*: = \Hom_k(L, k)$ with the $G$-action
defined as follows: for $g$ in  $G$, $\phi$ in $L^*$ and $x$ in $L$, 
$(g \phi) (x) := \phi(g^{-1}x)$.

\begin{prop} \label{prop:SWduality}
Let $G$ be a finite group. A map $f \colon M \rar N$ between $kG$-modules
is a dual ghost if and only if $f^* \colon N^* \rar M^*$ is a ghost.
\end{prop}

\begin{proof}
Consider the natural isomorphism
\[ \Hom_{kG} (N, \Hom_k(T, k)) \cong \Hom_{kG}(T, \Hom_k(N, k)), \]
where $T$ is a Tate resolution of $k$ and $N$ is any $kG$-module;
see~\cite[Proposition~3.1.8]{ben-1}, for instance. Since $\Hom_k(T, k)$ is a
complete injective resolution of $k$, taking homology of the chain complexes
in the last isomorphism gives natural isomorphisms
\[ \uHom_{kG}(N, \Omega^{-n} k) \cong \uHom_{kG} (\Omega^n k, N^*). \]
This implies that a map $d \colon M \rar N$ is a dual ghost if and only if
$d^* \colon N^* \rar M^*$ is a ghost.
\end{proof}

The second isomorphism used in the proof of the above proposition is
\emph{Spanier-Whitehead duality} for the stable module category.

We now use Tate duality to show  that one can use group cohomology ($\Tate^i(G,-)$, $i
\ge 0$) to detect ghosts.

\begin{thm} \label{prop:ghostbycohomology}
Let $G$ be a $p$-group. A map $f \colon M \rar N$ between finite-dimensional
$kG$-modules is a ghost if and only if the following two conditions hold:
\begin{enumerate}
\item $\Tate^i(G, f) \colon \Tate^i(G, M) \rar \Tate^i(G, N)$ is zero for all $i \ge 0$.
\item $\Tate^i(G, f^*) \colon \Tate^i(G, N^*) \rar \Tate^i(G, M^*)$ is zero for all $i \ge
0$.
\end{enumerate}
\end{thm}

\begin{proof}
Clearly it suffices to show that statement (2) is equivalent to the
statement:
\[\Tate^{-i-1}(G, f) \colon \Tate^{-i-1}(G, M) \rar \Tate^{-i-1}(G, N)\]
is zero for
all $i \ge 0$. Recall that Tate duality~\cite{Carlson-Tateduality} gives a natural isomorphism
\[\Tate^{-i -1} (G, L) \cong (\Tate^i (G, L^*))^*\]
for any finite-dimensional module $L$. Thus, for each $i \ge 0$, we have the
following commutative diagram, where the vertical maps are induced by $f$:
\[
\xymatrix{
 \Tate^{-i -1}(G, M) \ar[r]^{\cong} \ar[d] & (\Tate^i(G, M^*))^* \ar[d] \\
  \Tate^{-i -1}(G, N) \ar[r]^{\cong} & (\Tate^i(G, N^*))^*.
}
\]
 Since the two horizontal maps are isomorphisms, the left vertical map is
zero if and only if the right vertical map is zero. Finally, by the
faithfulness of the vector space duality functor, the right vertical map is zero
if and only if statement (2) holds. So we are done.
\end{proof}

Combining Spanier-Whitehead duality and Tate duality gives us the following
interesting result.

\begin{cor} \label{cor:ghost=dualghost}
A map $f \colon M \rar N$ between finite-dimensional $kG$-modules is a ghost
if and only if it is a dual ghost.
\end{cor}

\begin{proof}
By Proposition \ref{prop:SWduality}, we know that $f \colon M \rar N$ is a dual ghost
if and only if $f^* \colon N^* \rar M^*$ is a ghost. And by Proposition \ref{prop:ghostbycohomology},
$f^*$ is  a ghost if and only if
\begin{enumerate}
\item $\Tate^i(G, f^*) \colon \Tate^i(G, N^*) \rar \Tate^i(G, M^*)$ is zero for all $i \ge 0$, and
\item $\Tate^i(G, f) \colon \Tate^i(G, M) \rar \Tate^i(G, N)$ is zero for all $i \ge
0$.
\end{enumerate}
(We get the second statement from the fact that double dual $f^{**}$ is naturally isomorphic to $f$.)
The last two statements are in turn equivalent, again by Proposition \ref{prop:ghostbycohomology},
to the statement that $f$ is a ghost.
\end{proof}

\begin{rem}
The analogue of Corollary~\ref{cor:ghost=dualghost}  fails in the derived category of a
commutative ring. For example, in $D(\ints)$, let $X$ be the cone of the map
$\ints \stk{p} \ints$, and let $Y = \Sigma \ints$. Then the map $f \colon X
\rar Y$ which projects onto the top class, i.e., the map that is the identity in
degree $1$ and zero elsewhere, is easily seen to be a non-trivial ghost.
However, $f$ is  not a dual ghost.  In fact, the composite $X \stk{f} Y
\stk{=} Y (= \Sigma \ints)$ is just $f$, which is non-trivial.
\end{rem}

The first part of the following corollary gives an
alternative proof of Corollary~\ref{cor:GHwithtargetk-1}.

\begin{cor} \label{cor:GHwithtargetk-2}
Let $M$ be a finite dimensional $kG$-module. Then we have the following:
\begin{enumerate}
\item If $f \colon M \rar \Omega^i k$ is a ghost, then $f$ is stably trivial. In
other words, the GH with target $k$ holds.
\item If $f \colon \Omega^i k \rar M$ is a dual ghost, then $f$ is stably trivial.
\end{enumerate}
\end{cor}

\begin{proof}
If $f \colon M \rar \Omega^i k$ is a ghost, then by
Corollary~\ref{cor:ghost=dualghost}, $f$ is also a dual ghost, so the composite
$M \lstk{f} \Omega^i k \lstk{=} \Omega^i k$, which is just $f$, is stably
trivial.

The proof of the second statement is similar.
\end{proof}

\begin{rem}
Using the results in this section, we can show that ``the GH with domain $L$''
(the statement that every ghost in $\stmod(kG)$ with domain $L$ is trivial) holds
if and only if ``the GH with target $L^*$'' holds. This generalises
Corollary~\ref{cor:GHwithtargetk-2}. We leave the easy details to the reader.
\end{rem}

\section{Universal ghosts} \label{sec:universalghosts}

A ghost $\Phi \colon M \rar N$ between $kG$-modules is said to be a
\emph{universal ghost} if every ghost out of $M$ factors through $\Phi$. (Such
a map should technically be called \emph{weakly} universal, since we do not
assume the factorisation is unique.) We will show that there exists a
universal ghost out of any given $kG$-module. Let $M$ be a $kG$-module.
We assemble all the homogeneous elements in $\Tate^*(G, M)$ into a map
\[\bigoplus_{\eta\in\Tate^*(G, M)} \Omega^{|\eta|} k \lrar M , \]
where $|\eta|$ is the degree of $\eta$. Completing this map to an exact
triangle in $\StMod(kG)$, we get
\begin{equation} \label{eq:univ-ghost}
\bigoplus_{\eta\in\Tate^*(G, M)} \Omega^{|\eta|} k \lrar M \lstk{\Phi_M} U_M.
\end{equation}

We now recall a couple of easily established facts (see \cite{CCM2} for proofs) which we will need in the sequel.

\begin{prop}[\cite{CCM2}]\label{prop:univ}
The map $\Phi_M \colon M \rar U_M$ is a universal ghost out of $M$.
\end{prop}



Our next proposition characterises modules out of which all ghosts vanish.

\begin{prop}[\cite{CCM2}]\label{prop:trivial}
Let $M$ be a $kG$-module. Then the following are equivalent statements in the
stable module category:
\begin{enumerate}
\item All ghosts out of $M$ are trivial.
\item The universal ghost $\Phi_M\colon M \rar U_M$ is trivial.
\item $M$ is a retract of a direct sum of suspensions of the trivial representation.
\end{enumerate}
Moreover, if $M$ is finite-dimensional, (3) can be replaced with the
condition that $M$ is stably isomorphic to a finite direct sum of suspensions of the trivial
representation.
\end{prop}


Now suppose that $M$ is a finite-dimensional $kG$-module such that
$\Tate^*(G,M)$ is finitely generated as a graded module over $\Tate^*(G, k)$.
(This happens, for example, when $k$ is periodic, that is, $\Omega^i \, k$ is
stably isomorphic to $k$ for some $i \neq 0$.  See \cite{CarCheMin} for interesting
and non-trivial examples in the non-periodic case.) We will show that a universal
ghost out of $M$ can be constructed in the category $\stmod(kG)$; that is, the
target module of the universal ghost out of $M$ can be chosen to be
finite-dimensional as well. This is done as follows. Let $\{ v_j \}$ be a
finite set of homogeneous generators for $\Tate^*(G,M)$ as an
$\Tate^*(G,k)$-module. These generators can be assembled into a map
\[
\bigoplus_j \, \Omega^{|v_j|}\,k \lrar M
\]
in $\stmod(kG)$. This map can then be completed to a triangle
\[
\bigoplus_j \, \Omega^{|v_j|}\,k \lrar M \lstk{\Psi_M} F_M 
\]
in $\stmod(kG)$.
By construction, it is clear that the first map in the above triangle is
surjective on the functors $\uHom_{kG}(\Omega^l k, -)$ for each $l$. Therefore, the
second map $\Psi_M$ must be a ghost.
Universality of $\Psi_M$ is easy to see, so we have proved the following
proposition.

\begin{prop}\label{prop:fduniv}
Suppose $M$ is a finite-dimensional $kG$-module such that $\Tate^*(G,M)$
is finitely generated as a graded module over $\Tate^*(G, k)$.
Then a universal ghost out of $M$ can be constructed in $\stmod(kG)$.
In particular, this applies when $k$ is periodic.
\end{prop}

\begin{cor}\label{cor:finitelygenerated}
Let $M$ be a finite-dimensional module such that $\Tate^*(G,M)$ is finitely
generated as a graded module over $\Tate^*(G, k)$. Then the following are
equivalent statements in the stable module category:
\begin{enumerate}
\item All ghosts out of $M$ are trivial.
\item All ghosts out of $M$ into finite-dimensional modules are trivial.
\item The universal ghost $\Psi_M\colon M \rar F_M$ is trivial.
\item $M$ is  a finite direct sum of suspensions of the trivial representation.
\end{enumerate}
\end{cor}


We now give some applications of universal ghosts.
We begin with a characterisation of finite-dimensional indecomposable projective-free representations that are
isomorphic to a Heller shift of the trivial representation.

\begin{cor}\label{cor:ghostcondn} Let $G$ be a finite $p$-group and
let $M$ be a finite-dimensional indecomposable projective-free $kG$-module.
Then all ghosts out of $M$ are trivial if and only if $M \cong \wOmega^i\, k$
for some integer $i$.
\end{cor}

\begin{proof}
This follows from Proposition~\ref{prop:trivial} and the Krull-Schmidt theorem.
\end{proof}

\begin{cor}\label{cor:modulo}
Let $G$ be a finite $p$-group and let $M$ be a finite-dimensional
indecomposable projective-free $kG$-module. If $\dim(M)$ is not congruent to
$+1$ or $-1$ modulo $|G|$, then there exists a non-trivial ghost out of
$M$.
\end{cor}

\begin{proof}
By the previous corollary it suffices to show that $M \ncong \wOmega^i\,k$
for any $i$. This will be shown by proving that under the given hypothesis the
dimensions of the Heller shifts $\wOmega^i\, k$ are congruent to $+1$ or $-1$
modulo $|G|$.
Recall that $\wOmega^1\, k$ is defined to be the kernel of the augmentation
map, so we have a short exact sequence
\[ 0 \lrar \wOmega^1\, k \lrar kG \lrar k \lrar 0, \]
which tells us that $\dim (\wOmega^1\, k) \equiv  -1$ modulo $|G|$.
Inductively, it is clear from the short exact sequences
\[ 0 \lrar \wOmega^{i+1}\,k \lrar (kG)^t \lrar \wOmega^i\, k \lrar 0\]
that $\dim (\wOmega^i\,k) \equiv  (-1)^i$ modulo $|G|$ for $i \geq 0$. (Here
$(kG)^t$, for some $t$, is a minimal projective cover of $\wOmega^i \, k$.)
Also, since $\wOmega^i\,k \cong (\wOmega^{-i}\,k)^*$ in $\Mod(kG)$, it follows
that $\dim (\wOmega^i\,k) \equiv  (-1)^i$ modulo $|G|$ for each integer $i$.
In particular, $M \ncong \wOmega^i\, k$ for any integer $i$.
\end{proof}

\section{The ghost projective class} \label{sec:projectiveclass}

In order to measure the degree to which the GH fails, it is natural to
consider the smallest integer $l$ such that the composite of any $l$ ghosts
between finite-dimensional $kG$-modules is trivial in the stable module category. 
We will show (see Theorem~\ref{thm:3inequalitiesofnumbers}) 
that such an integer always exists. This integer will be called the \emph{ghost number of $\stmod(kG)$}, or, more briefly, the
\emph{ghost number of $kG$} and can be best understood using the
concept of a projective class. So we begin with a recollection of the notion
of a projective class in a triangulated category. A good reference for this is
\cite{ch}, where projective classes were studied in the stable homotopy
category and the derived category of a ring.

\subsection{Projective classes}
Let $\T $ denote a triangulated category. A \emph{projective class} in $\T$ is a pair $(\mcP , \G )$ where $\mcP$
is a class of objects and $\G$ is a class of maps in $\T$ which satisfy the following properties:
\begin{enumerate}
\item The class of all maps $X \rar  Y$ such that the composite $P \rar X \rar Y$ is zero for all $P$ in $\mcP$ and all maps
$P \rar X$ is precisely $\G$.
\item The class of all objects $P$ such that the composite $P \rar X \rar Y$ is zero for all maps
$X \rar Y$ in $\G$ and all maps
$P \rar X$ is precisely $\mcP$.
\item For each object $X$ there is an exact triangle $P \rar X \rar Y$ with  $P$ in $\mcP$ and $X \rar Y$ in $\G $.
\end{enumerate}

It follows that the maps in $\G $ form an ideal in $\T$.
That is, if $f$ and $g$ are parallel maps in $\G $, then $f+g$ is in $\G $, and if $f$, $g$, and $h$ are composable
with $g$ in $\G $, then both $fg$ and $gh$ are also in $\G $.

Once we have a projective class as defined above, we can form ``derived'' projective classes in a natural way as follows.
The powers $\G^{n}$ of the ideal $\G $ form a decreasing filtration of the maps in $\T $, and
each  $\G^{n}$ is part of a projective class.  The corresponding classes of objects are obtained as follows.
Let $\mcP^{1} = \mcP$ and inductively define $\mcP^{n}$ to be the collection of retracts of objects $M$ that appear
in a triangle
\[
 A \lrar M \lrar B,
\]
where $A$ belongs to $\mcP^{1}$ and $B$ belong to $\mcP^{n-1}$. The classes $\mcP^{n}$ form an increasing filtration
of the objects in $\T $. It is a fact \cite[Theorem 1.1]{ch} that $(\mcP^{n}, \G^{n})$ is a projective class for each $n$.
We set $\mcP^{0}$  to be the collection of zero objects in $\T $ and $\G^{0}$ to be the collection of all maps in $\T $.
$(\mcP^0, \G^0)$ is the \emph{trivial projective class}.

\subsection{The ghost projective class}
Now we specialise to the stable module category to define the \emph{ghost
projective class} in $\StMod(kG)$.  
The ideal $\G$ consists of the class of ghosts.
The associated class $\mcP$ of objects consists of retracts of direct sums of
suspensions of the trivial representation.

\begin{prop}
The pair $(\mathcal{P}, \mathcal{G})$ is a projective class in $\StMod(kG)$.
\end{prop}

\begin{proof}
It is clear that $\mathcal{P}$ and $\mathcal{G}$ are orthogonal, i.e., the
composite $P \rar M \stk{h} N$ is zero for all $P$ in $\mathcal{P}$,  for all
$h$ in $\mathcal{G}$, and all maps $P \rar M$. So by  \cite[Lemma 3.2]{ch} it
remains to show that for all $kG$-modules $M$, there exists a triangle $P \rar
M \rar N$ such that $P$ is in $\mathcal{P}$ and $M \rar N$ is in
$\mathcal{G}$.  The universal ghost \eqref{eq:univ-ghost} out of $M$ has this
property, so we are done.
\end{proof}

In some special cases one can also build a ghost projective class in $\stmod(kG)$. Let $\mcP_c$ denote the
collection of finite direct sums of suspensions of $k$  and $\G_c$ the class of ghosts in $\stmod(kG)$.
(Note that the collection $\mcP_c$ is already closed under retractions by the Krull-Schmidt theorem.) Then we have the
following proposition.

\begin{prop}\label{prop:compact-projective-class}
Let $G$ be a finite $p$-group such that the trivial module $k$ is periodic.
Then $(\mcP_c, \G_c)$ is a projective class in $\stmod(kG)$.
\end{prop}

The proof below only requires that each finite-dimensional module
has finitely generated Tate cohomology.  
It is shown in~\cite{CarCheMin} that this is equivalent to $G$ being 
cyclic or a generalised quaternion group, and by a result of 
Artin and Tate~\cite[p.\ 262]{CarEil} this is equivalent to $k$
being periodic.

\begin{proof}
Orthogonality of $\mcP_c$ and $\G_c$ is clear. We have already seen in Proposition~\ref{prop:fduniv} that under the
given hypothesis a universal ghost out of a finite-dimensional module can be constructed within $\stmod(kG)$.
So the proposition follows from \cite[Lemma 3.2]{ch}.
\end{proof}


Whether or not the pair $(\mcP_c, \G_c)$  forms a projective class, we can define $\mcP_c^m$ and $\G_c^m$ as in the 
previous subsection.  A finite-dimensional $kG$-module is said to have
\emph{generating length}
$m$ if it belongs to $\mcP^m_{c}$ but not to $\mcP^{m-1}_c$, and
\emph{ghost length} $m$ if it is the domain of a non-zero map in $\G_c^{m-1}$ but
not in $\G_c^{m}$.

We now prove a sequence of inequalities.

\begin{prop} \label{prop:generating>=ghost}
Let $G$ be a finite $p$-group and $M$ a finite-dimensional $kG$-module.
Then
\[ \text{ghost length of } M \le \text{generating length of } M .\]
Moreover, equality holds if $(\mcP_c, \G_c)$ is a projective class.
\end{prop}

\begin{proof}
If the generating length of $M$ is one, then $M$ is a finite direct
sum of suspensions of $k$. Then clearly all ghosts out of $M$ vanish, which means the ghost
length of $M$ is also one.  Suppose that the generating length of $M$ is two. Then $M$ can be
chosen to be a direct summand of $M'$ which can be obtained as an extension of finite dimensional modules
\[ 0 \lrar  \oplus_I \Omega^i k \lstk{\alpha}  M' \lstk{\beta} \oplus_J \Omega^i k \lrar 0.  \]
It suffices to show that every two fold composite
\[M' \lstk{f} A \lstk{g} B \]
out of $M'$ is trivial in the stable category, for then the same would be true
for $M$. Consider the following commutative diagram in $\stmod(kG)$:
\[
\xymatrix{
\oplus_I \Omega^i k \ar[r]^{\alpha} \ar[dr]_{f \alpha = 0} & M' \ar[r] ^{\beta} \ar[d]^f &  \oplus_J  \Omega^i k
\ar@{.>}[dl]_{\hspace{7 mm}\tilde{f}} \ar[ddl]^{g \tilde{f} = 0} \  \\
& A \ar[d]^g &   \\
& B & .
}
\]
Since $f$ is a ghost, $f \alpha = 0$, therefore the map $f$ factors as $f = \tilde{f} \beta$. But then
\[ g f = g (\tilde{f} \beta) = (g \tilde{f}) \beta = 0 \beta = 0.\]
(The third equality follows because $g$ is a ghost.)
Since $gf=0$, the ghost length of $M$ is at most two.
The induction step is similar. 

The last statement of the proposition follows
directly from~\cite[Prop.~3.3]{ch}.
\end{proof}

\begin{prop} \label{prop:inequalityoflengths}
Let $G$ be a finite $p$-group and $M$ a finite-dimensional $kG$-module.
Then
\[  \text{generating length of } M
\le \text{radical length of } M. \]
\end{prop}

\begin{proof}
The Jacobson radical $J=J(kG)$ of $kG$ is nilpotent. The radical length of
$M$ is the smallest integer $h$ such that $J^h M = 0$.
This gives the radical or Lowey series for $M$:
\[ M \varsupsetneq JM \varsupsetneq  \cdots \varsupsetneq J^{h-1} M \varsupsetneq J^h M = 0.\]
Note that $J$ annihilates each successive quotient and hence each of them is
a direct sum of trivial representations. This shows that the
generating length is at most the radical length.
%
\end{proof}

\begin{rem}
Propositions~\ref{prop:generating>=ghost} and~\ref{prop:inequalityoflengths}
show that if the radical length of $M$ is $l$, then any map $M \rar N$
which is a composite of $l$ ghosts is trivial.
In contrast, it follows from Corollary~\ref{cor:soc-rad} that any map $L \rar M$
which is a composite of $l$ ghosts is trivial.
\end{rem}

Recall that the \emph{nilpotency index} of the Jacobson radical $J(kG)$ of $kG$
is the smallest integer $m$ such that $J(kG)^m = 0$.

\begin{prop} \label{prop:nilpotency>radical}
Let $G$ be a finite $p$-group and $M$ a
projective-free $kG$-module. Then
\[  \text{radical length of } M < \text{nilpotency index of } J(kG) \le |G|. \]
\end{prop}

\begin{proof} Let $m$ be the nilpotency index of $J(kG)$.
We begin by noting that since $G$ is a $p$-group, the last non-zero power $J(kG)^{m-1}$
of $kG$ is  the unique non-zero minimal ideal in $kG$; see \cite[p.\ 92]{ben-1}.
For the first inequality, it is enough to show that $J(kG)^{m-1}M = 0$.
Let $x$ be an element of $M$.  Since $M$ is projective-free,
$\Ann(x) \ne 0$.  Thus $\Ann(x)$
contains the unique non-zero minimal ideal $J(kG)^{m-1}$ of $kG$.
That is, $J(kG)^{m-1} x = 0$.

Since the powers of $J(kG)$ form a strictly decreasing series, then
the nilpotency index of $J(kG)$ is at most $\dim_k kG = |G|$.
\end{proof}

These inequalities show that for each $p$-group $G$, the ghost lengths and
generating lengths
of finite-dimensional $kG$-modules are uniformly bounded above.
So we define the \emph{generating number} of $kG$ to be the least upper
bound of the generating lengths of all finite-dimensional $kG$-modules,
and the \emph{ghost number} of $kG$ to be the least upper bound of
the ghost lengths of all finite-dimensional $kG$-modules.
We don't know if the generating number and ghost number of $kG$ depend
on the field $k$.

Combining the above results gives:

\begin{thm}\label{thm:3inequalitiesofnumbers}
Let $G$ be a finite $p$-group. Then
\begin{align*}
\text{ghost number of } kG  & \le \text{generating number of } kG  \\
                                    &  < \text{nilpotency index of } J(kG) \le |G|.
\end{align*}
In particular, the generating number and ghost number of the group
algebra of any finite $p$-group are finite, and any
composite of $|G|-1$ ghosts in $\stmod(kG)$ is trivial.
\end{thm}

\begin{rem}
A similar argument involving the projective class $(\mcP, \G)$ shows
that any composite of $m - 1$ ghosts in $\StMod(kG)$ is trivial, where
$m$ is the nilpotency index of $J(kG)$.
\end{rem}

\begin{prop} \label{prop:liftingghostnumber}
Let $H$ be a subgroup of a finite $p$-group $G$. Then
\[ \text{ghost number of } kH \le \text{ghost number of }kG. \]
\end{prop}

\begin{proof}
In \cite{CCM3}, we have shown that the induction functor
\[ \uInd \colon \stmod(kH) \lrar \stmod(kG), \]
which sends a $kH$-module $M$ to $M \up := kG \otimes_{kH} M$, preserves ghosts
and non-trivial maps. It follows that the ghost number of $kH$ is no more
 than that of $kG$.
\end{proof}

\section{Computing ghost numbers} \label{sec:computing}

 We now investigate the problem of computing ghost numbers and generating
numbers of some specific groups. The following lemma that we learned from Dave Benson~\cite{CCM3}
will be very helpful in computations.

If an element $\theta$ belongs to the centre of $kG$, there is a natural self $kG$-linear map on any $kG$-module $M$
given by left multiplication by $\theta$. We will denote this map by $\theta \colon M \rar M$.

\begin{prop} \label{prop:conjugacyghost}
Let $G$ be a finite $p$-group and let $M$ be a $kG$-module. If an element $\theta$
belonging to $J(kG)$ is central in $kG$, then the map
\[ \theta \colon M \lrar M \]
is a ghost.
In particular, if $g \in G$ is central, then the map
\[ g-1 \colon M \lrar M \]
is a ghost.
\end{prop}

\begin{proof}
The proof of~\cite[Lemma~2.2]{CCM3} applies without change.
\end{proof}


\subsection{Cyclic $p$-groups}
Recall that for cyclic groups, the trivial representation is periodic, and
so by Proposition~\ref{prop:compact-projective-class}
$(\mcP_c, \G_c)$ forms a projective class. So the ghost length and generating length are the same for
modules over cyclic $p$-groups.

\begin{prop}\label{ghostbound-cyclic}
All finite-dimensional $k\,C_{p^r}$-modules have ghost length at most $\lceil(p^{r}-1)/2\rceil $.
\end{prop}

Here $\lceil y \rceil$ denotes the smallest integer that is greater than or equal to $y$.

\begin{proof}
Since the characteristic of $k$ is $p$, we have $kC_{p^r} \cong k[x]/(x^{p^r})$,
with $x$ corresponding to $\sigma - 1$, where $\sigma$ is a generator of $C_{p^r}$.
A finite-dimensional indecomposable projective-free module over
$k[x]/(x^{p^r})$ is of the form $k[x]/x^i$ for
$1 \le i \le p^r-1$. It is also clear that $\wOmega(k[x]/(x^i)) \cong k[x]/(x^{p^r-i})$. This tells us that
the ghost length of $k[x]/(x^i)$ is the same as that of $k[x]/(x^{p^r-i})$. For $1 \le i
\le \lceil(p^{r}-1)/2 \rceil$, we show that the ghost length of $k[x]/x^i$ is at most $i$. For this
it is enough to observe that we have short exact sequences
\[ 0 \lrar k \lrar k[x]/x^i \lrar k[x]/x^{i-1} \lrar 0\]
of modules over $k[x]/x^{p^r}$,
for $2 \le i \le \lceil(p^{r}-1)/2 \rceil$.
\end{proof}

\begin{prop}
There exists a composable sequence of $\lceil(p^{r}-1)/2\rceil-1$ ghosts in\break $\stmod (kC_{p^r})$
whose composite is non-trivial.
\end{prop}

\begin{proof}
Recall that $kC_{p^r} \cong k[x]/(x^{p^r})$.
Let $h \colon k[x]/x^d \rar k[x]/x^d$ be multiplication by $x = \sigma - 1$,
where $d = \lceil(p^{r}-1)/2\rceil$.
By Proposition~\ref{prop:conjugacyghost}, $h$ is a ghost.
To see that $h^{d - 1}$ is non-trivial, we have to show that it cannot factor
through the projective cover
$k[x]/x^{p^r} \twoheadrightarrow k[x]/x^d$, i.e., that we cannot have a commutative diagram
\[
\xymatrix{
k[x]/x^d \ar[rr]^{x^{d - 1}} \ar@{-->}[dr]& & k[x]/x^d  \\
& k[x]/x^{p^r} \ar@{>>}[ur] & .
}
\]
By considering the images of the generator of the left-hand cyclic module in the above diagram, one can
easily see that the existence of such a factorisation would mean that
\[
(d-1) + (d-1) \ge p^r-1 ,
\]
or, equivalently, that
\[\lceil(p^{r}-1)/2\rceil-1 + \lceil(p^{r}-1)/2\rceil-1 \ge p^r -1 . \]
It is straightforward to verify that this inequality fails for all primes $p$ and all positive integers $r$.
So we are done.
\end{proof}

Combining these two propositions, we get the following theorem.

\begin{thm}\label{ghostnumber-cyclic}
The ghost number of $kC_{p^{r}}$ is $\lceil(p^{r}-1)/2\rceil$.
\end{thm}

\begin{cor}[\cite{CCM3}]
The GH holds for $kC_{p^{r}}$ if and only if $p^{r}$ is equal to  $2$ or $3$.
\end{cor}

\begin{proof}
Recall that the GH holds for $kG$ precisely when  the
ghost number of $kG$ is  $1$. So, by Theorem~\ref{ghostnumber-cyclic}, we conclude that the GH holds for $kC_{p^r}$
if and only if $\lceil(p^{r}-1)/2\rceil =  1$. The last equation holds
if and only if  $p^r = 2$ or $3$.
\end{proof}

\subsection{The Klein four group}

\begin{prop} \label{prop:V_4genlength}
Let $M$ be a finite-dimensional indecomposable projective-free $kV_4$-module.
Then we have the following:
\begin{enumerate}
\item  If $M$ is odd-dimensional, then it has generating length one.
\item  If $M$ is even-dimensional, then it has generating length two.
\end{enumerate}
\end{prop}

\begin{proof}
It is well-known that the odd-dimensional indecomposable modules are precisely
the Heller shifts of the trivial representation; see~\cite[Theorem 4.3.2]{ben-1}, for
instance. So they all have generating length one by definition.
Since $M$ is projective-free, one can show using the
classification of the indecomposable $kV_4$-modules (e.g.,~\cite[Theorem 4.3.2]{ben-1}), or directly,
that there is a short exact sequence
\[ 0 \lrar  M^{V_4}  \lrar M \lrar M_{V_4}   \lrar 0,\]
where the invariant submodule $M^{V_4}$ and the coinvariant module $M_{V_4}$
are both direct sums of the trivial representation $k$.  Thus $M$ has 
generating length at most two.  Moreover, if $M$ is even-dimensional and
indecomposable, then  $M$  is not isomorphic to $\wOmega^i k$ for any $i$. In
particular, $M$ cannot have generating length one. So we are done.
\end{proof}


\begin{thm} \label{thm:compositionofghostsV_4}
The ghost number and the generating number of $kV_4$ are both two.
\end{thm}

\begin{proof}
Since every finite-dimensional module is a sum of indecomposables, the
statement about the generating number follows from the above proposition.
Since the ghost number is at most the generating number, we
only have to show that the ghost number is bigger than one. This
follows from~\cite{CCM3} because there we showed that the GH  fails for
$\stmod(kV_4)$.
\end{proof}

\subsection{The quaternion group}

By our main result on the GH in~\cite{CCM3}, we know that the GH fails for the
quaternion group $Q_8$ of order $8$. Now we give bounds on the ghost number.
Since the trivial representation of $Q_8$ is periodic, we know from
Proposition~\ref{prop:compact-projective-class} that the ghost
projective class exists in $\stmod(kQ_8)$. Therefore the ghost number and the
generating number of $kQ_8$ are the same.

\begin{prop}
The ghost number of $k Q_8$ is at least two and at most four.
\end{prop}

\begin{proof}
The nilpotency index of $J(kQ_8)$ is 5 (see Subsection~\ref{subse:nilp-index}), so
by Theorem~\ref{thm:3inequalitiesofnumbers}, we know
that the ghost number of $kQ_8$ is at most $4$. We have already seen that the GH fails for $Q_8$, so
the ghost number is at least 2.
\end{proof}

In the following example we will give another disproof of the GH for the group
$Q_8$ by exhibiting an explicit finite-dimensional module with
ghost length two. This example should also illustrate some of the ideas surrounding
projective classes.

\begin{example} \label{counterexample-GH-Q8}
A minimal presentation for $Q_8$ is given by
\[ Q_8 = \la x, y \; |\; x^4 = 1, x^2 = y^2, yxy^{-1} = x^{-1}  \ra . \]
The structure of the left $kQ_8$-module $J(kQ_8)^3$ can be obtained using
Jennings' theorem~\cite{jennings} or otherwise. It is shown in the diagram below:
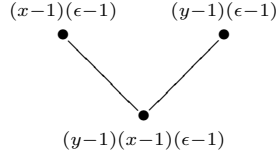
\[
\xymatrix{
*{\bullet} \save[]+<0pt,9pt>*{\scriptstyle (x - 1)(\epsilon -1)}\restore \ar@{-}[dr]
  & & *{\bullet} \save[]+<0pt,9pt>*{\scriptstyle (y - 1)(\epsilon -1)}\restore \ar@{-}[dl] \\
& *{\bullet} \save[]+<0pt,-9pt>*{\scriptstyle (y - 1)( x - 1)(\epsilon -1)}\restore &
}
\]
Here a bullet is a one-dimensional $k$-vector space,
the southwest line segment corresponds to the action of $x - 1$,
the southeast line segment corresponds to the action of $y -1$,
and if no line segment emanates from a bullet in a given
direction, then the corresponding action is trivial.
$\epsilon$ is the central element $x^2 (= y^2)$.

It is clear from the diagram that the invariant submodule of $J(kQ_8)^3$ is
one-dimen\-sional and therefore we conclude that $J(kQ_8)^3$ is
indecomposable (see, for example,~\cite[Lemma~3.2]{CCM2}).
Also note that $J(kQ_8)^3$ is projective-free. Moreover, the
dimension of $J(kQ_8)^3$ is 3, which is neither $+1$ nor $-1$ modulo 8
($=|Q_8|$). Thus by Corollary~\ref{cor:modulo}  we know that there exists a
non-trivial ghost in $\StMod(kQ_8)$ whose domain is $J(kQ_8)^3$. 
By Corollary~\ref{cor:finitelygenerated}, the target can be chosen
to lie in $\stmod(kQ_8)$.
In particular, the ghost length of $J(kQ_8)^3$ is at least 2. On the
other hand the generating length of $J(kQ_8)^3$ is at most 2 because
we have a short exact sequence of $kQ_8$-modules
\[ 0 \lrar  k \lrar J(kQ_8)^3 \lrar k \oplus k \lrar 0 . \]
Since the ghost length is always less than or equal to the
generating length, we conclude that
the ghost length and the generating length of $J(kQ_8)^3$ are both 2.
\end{example}

\subsection{Ghost numbers of abelian groups}

We begin with an extremely useful proposition.
This is a  slight generalisation of a result we learned from
Dave Benson that appeared in~\cite{CCM3}.

\begin{prop} \label{prop:GHranktwogroups}
Let $G$ be a finite $p$-group and let $H$ be a non-trivial proper subgroup
of $G$.  Let $\theta = \sum_g \alpha_g g$ be a central element of $kG$
such that $\sum_{h \in H} \alpha_h \neq 0$.
Then multiplication by $\theta$ on $k_H \up$
is stably non-trivial, where $k_H$ is the trivial $kH$-module.
In particular, if $g$ is a central element in $G - H$,
then multiplication by $g - 1$ on $k_H \up$ is a non-trivial ghost.
\end{prop}

We include a proof, since this is slightly more general
than~\cite[Lemma~2.3]{CCM3}:
we do not require $H$ to be normal and
we include more general $\theta$.

\begin{proof}
Recall that $k_H \up$ denotes the induced module $kG \otimes_{kH} k_H$, and
that induction is both left and right adjoint to restriction. These
adjunctions give rise to natural $kH$-linear maps $k_H \rar k_H \up \down$, sending $x$ to
$1 \otimes x$, and $k_H \up \down \rar k_H$, sending $g \otimes x$ to $x$ if
$g \in H$ and to $0$ otherwise.

To show that $\theta \colon k_H\up \rar k_H\up$ is stably non-trivial,
it is enough to show that
\[\theta\down \colon k_H\up\down \rar k_H\up\down\]
is stably non-trivial.  For this, it is enough to show that the
composite
\[k_H \rar k_H\up\down \lstk{\theta\down} k_H\up\down \rar k_H\]
is stably non-trivial.  But this composite is multiplication
by $\sum_{h \in H} \alpha_h$, which is non-zero by assumption.
And since $H$ is non-trivial, all non-zero maps $k_H \rar k_H$
are stably non-trivial.

The last statement follows from the first part of this proposition,
combined with Proposition~\ref{prop:conjugacyghost}.
\end{proof}

\begin{thm} \label{thm:boundsforabelian}
Let $G$ be an abelian $p$-group and let $m$ denote the nilpotency index of
$J(kG)$. Then we have
\[ m - p^{r-1}(p-1) \le \text{ghost number of } kG \le m-1 , \]
where $p^{r}$ is the order of the smallest cyclic summand of $G$.
\end{thm}

\begin{proof}
We have already seen the upper bound, so we only have to establish the lower
bound. Let $C_{p^r}$ be the smallest cyclic summand of $G$, so that for some
integers $r_i \ge r$,
\[ G =  C_{p^r} \oplus C_{p^{r_1}} \oplus \cdots \oplus C_{p^{r_t}}.\]
Let $H$ be the subgroup of order $p$ in $C_{p^r}$ and let $k_{H}$ be the trivial
$kH$-module. Set $M = k_H \up$. We will produce $m - p^{r -1}(p-1) - 1$
ghosts $M \rar M$ whose composite is stably non-trivial. This will give the desired
lower bound for the ghost number. Let $g$ be a generator for $C_{p^r}$
and let $g_i$ be a generator for $C_{p^{r_i}}$ for each $i$.
By Proposition~\ref{prop:conjugacyghost}, the map
\[ \theta = (g-1)^{p^{r-1}-1} (g_1-1)^{p^{r_1}-1} \cdots (g_t-1)^{p^{r_t}-1}
   \colon M \rar M \]
is a composite of
\[ (p^{r-1} - 1) + (p^{r_1} -1) + \cdots + (p^{r_t} - 1)\]
ghosts.
One can see easily that the nilpotency index of $J(kG)$ is
\[ m = 1 + (p^{r} - 1) + (p^{r_1} - 1 ) + \cdots + (p^{r_t} - 1).\]
Thus $\theta$ is a composite of $m - p^{r -1}(p-1) - 1$ ghosts.
Now to see that $\theta$ is stably non-trivial on $M$, it is enough to
note that if $\theta \in kG$ is written $\sum_{g \in G} \alpha_g g$, then
$\alpha_h = 0$ for $h \in H$ unless $h = e$.  So
Proposition~\ref{prop:GHranktwogroups} applies.
\end{proof}


We derive some easy corollaries.

\begin{cor} \label{cor:C_2splitsummand}
Let $G$ be an abelian $2$-group which has $C_2$ as a summand. Then the ghost
number of $kG$ is one less than the nilpotency index of $J(kG)$.
\end{cor}

\begin{proof}
In this case, both the lower bound and the upper bound for the 
ghost number of $kG$ are one less than the nilpotency index
of $J(kG)$.
\end{proof}

\begin{cor} \label{cor:elementaryabelian2group}
Let $G$ be an elementary abelian $2$-group of rank $l$, i.e., $G \cong
(C_2)^l$. Then the ghost number of $kG$ is $l$.
\end{cor}

\begin{proof}
The nilpotency index of $J(k(C_2)^l)$ is easily shown to be $l + 1$.
\end{proof}

\begin{thm}
Let $G$ be an abelian $p$-group. The  ghost number of $kG$ is $2$ if and only if
$G$ is $C_4$, $C_2 \oplus C_2$, or $C_5$.
\end{thm}

\begin{proof}
By Theorem~\ref{ghostnumber-cyclic} and Corollary~\ref{cor:elementaryabelian2group}
we know that the three given groups have
ghost number $2$. Now we prove the converse. An easy exercise using the
structure theorem  tells us that if $|G| > 5$, then $G$ contains one of the
following groups as a subgroup: $C_{p^l}$ $(p^l > 5)$, $C_2 \oplus C_2 \oplus
C_2$,  $ C_p \oplus C_p$ $(p > 2)$, or $C_2 \oplus C_4$.  It is easily seen
using the lower bound in Theorem~\ref{thm:boundsforabelian} that
the ghost number of each of the above groups is at least $3$. Therefore, by
Proposition~\ref{prop:liftingghostnumber}, the ghost number of $G$ is also at
least $3$. So if the ghost number of $G$ is at most $2$, then $|G|$ should be
at most $5$. We know that $C_2$ and $C_3$ have ghost number $1$, and
the only remaining groups of order at most $5$ are
$C_4$, $C_5$ and $C_2 \oplus C_2$.
\end{proof}

\subsection{Ghost numbers for non-abelian groups}\label{subse:nilp-index}
Computing ghost numbers for non-abelian groups is much harder. However, the nilpotency indices
of $J(kG)$ can be computed using a formula of Jennings~\cite{jennings}.  For example, if
$G$ is a non-cyclic group of order $2^n$ which has a cyclic subgroup
of index $2$, then the nilpotency index of $J(kG)$ is $2^{n-1}+1$; see \cite{erdmann} for instance. Therefore
the composite of any $2^{n-1}$ ghosts in $\stmod(kG)$ is trivial.
By a well-known theorem~\cite[Ch.\ 8, pp.\ 134--135]{burnside}  every non-abelian $2$-group that
has a cyclic subgroup of index $2$ is a dihedral, semidihedral,
modular or  quaternion group. Thus we have upper bounds for the ghost numbers of their group algebras.

\end{document}